\documentclass{srthesis}
\pdfoutput=1

\usepackage[backend = biber, style = alphabetic]{biblatex}
\DeclareMathOperator*{\colim}{colim}

\newcommand*{\bcoprod}{+}

\newcommand*{\ul}{\underline}

\newcommand*{\sep}[1]{\bar{#1}\vphantom{#1}}

\newcommand*{\earg}{\,\cdot\,}
\newcommand*{\blt}{\bullet}

\newcommand*{\simp}{\mathrm{s}}
\newcommand*{\Set}{\mathrm{Set}}
\newcommand*{\sSet}{\simp\Set}

\newcommand*{\Topo}{\mathrm{Top}}
\newcommand*{\sTopo}{\simp\Topo}
\newcommand*{\Prof}{\mathrm{Prof}}
\newcommand*{\sProf}{\simp\Prof}
\newcommand*{\Extr}{\mathrm{Extr}}

\newcommand*{\TotDisc}{\mathrm{TotDisc}}

\newcommand*{\comp}{\uppi}

\newcommand*{\Hy}{\mathrm{Hy}}

\newcommand*{\Hywc}{\mathrm{Hy_{wc}}}
\newcommand*{\Hyawc}{\mathrm{Hy_{awc}}}
\newcommand*{\Ho}{\mathrm{Ho}}

\newcommand*{\et}{\mathrm{\acute{e}t}}
\newcommand*{\pe}{\mathrm{p\acute{e}}}
\newcommand*{\aff}{\mathrm{aff}}

\newcommand*{\Sh}{\mathrm{Sh}}
\newcommand*{\AbSh}{\mathrm{AbSh}}

\newcommand*{\NDeg}{\mathrm{N}}

\newcommand*{\rep}{\mathrm{h}}

\newcommand*{\ZZ}{\boldsymbol{\mathrm{Z}}}

\newcommand*{\Gal}{\mathrm{Gal}}

\newcommand*{\psh}{\mathrm{p}}
\newcommand*{\sh}{\mathrm{s}}
\newcommand*{\op}{\mathrm{op}}

\newcommand*{\Homol}{\mathrm{H}}

\newcommand*{\Spec}{\mathrm{Spec}}
\newcommand*{\const}{\mathrm{c}}
\newcommand*{\Hom}{\mathrm{Hom}}

\newcommand*{\sk}{\mathrm{sk}}
\newcommand*{\cosk}{\mathrm{cosk}}

\newcommand*{\Pro}{\mathrm{Pro}}

\newcommand*{\CoEq}{\mathrm{CoEq}}

\newcommand*{\LDer}{\mathrm{L}}
\newcommand*{\awc}{\mathrm{awc}}
\newcommand*{\Cond}{\mathrm{Cond}}
\newcommand*{\sCond}{\mathrm{sCond}}
\newcommand*{\How}{\mathrm{Ho_w}}
\newcommand*{\fw}{\mathrm{fw}}

\newcommand*{\sA}{\mathscr{A}}
\newcommand*{\sC}{\mathscr{C}}
\newcommand*{\sD}{\mathscr{D}}
\newcommand*{\sF}{\mathscr{F}}

\newcommand*{\sO}{\mathscr{O}}
\newcommand*{\sP}{\mathscr{P}}

\newcommand*{\mfi}{\mathfrak{i}}
\newcommand*{\mfp}{\mathfrak{p}}
\newcommand*{\mfS}{\mathfrak{S}}

\newcommand{\ulsi}[1]{\!\underline{\,{#1}}}

\newcommand*{\hsmash}{\mathrlap}

\newcommand*{\rto}{\rightarrow}
\newcommand*{\rrto}{\rightrightarrows}

\newcommand*{\rlto}{\rightleftarrows}

\begin{document}

\begin{center}

\LARGE
\textbf{The Pro-\'Etale Homotopy Type}

\vspace*{0.5cm}

\normalsize
Paul Meffle

\end{center}

\vspace*{7ex}

\renewcommand{\abstractname}{Abstract}

\begin{abstract}

In this paper we define the pro-\'etale homotopy type of a scheme and prove some of its expected properties.
Our definition is similar to the definition of the \'etale homotopy type by Michael Artin and Barry Mazur.
We prove that for a qcqs scheme the pro-\'etale homotopy type is profinite, determined by a single split affine wc hypercovering
and computes the cohomology of a certain class of sheaves.
We show that the pro-\'etale homotopy type of a w-contractible scheme is trivial and compute the pro-\'etale homotopy type
of the real numbers.
Moreover, we prove that a suitable version of $\comp_0$ composed with the pro-\'etale homotopy type gives back the space
of components of the base scheme.
We make some progress towards describing the pro-\'etale homotopy type of arbitrary fields.
Lastly, we give a refined definition of the pro-\'etale homotopy type using the theory by Ilan Barnea and Tomer M.~Schlank
and the theory of condensed sets by Dustin Clausen and Peter Scholze.
This allows us to define pro-\'etale homotopy groups associated to pointed qcqs schemes.

\end{abstract}

\setlength{\cftbeforetoctitleskip}{3.5ex plus 1ex minus .2ex}
\setlength{\cftaftertoctitleskip}{2.3ex plus .2ex}
\renewcommand{\cfttoctitlefont}{\large\bfseries}
\setlength{\cftbeforechapskip}{0pt}
\renewcommand{\cftchapfont}{\normalfont}
\renewcommand{\cftchappagefont}{\normalfont}

\tableofcontents

\clearpage

\chapter{Introduction}

The definition of the \'etale homotopy type in \cite{am69} by Michael Artin and Barry Mazur relies on the fact that the
\'etale site is locally connected.
As we show in Example~(\ref{exa:pe-morph-from-non-lc}) the pro-\'etale site is not locally connected.
As a result, we cannot apply the theory from \cite{am69} directly.
To address the issue we use the \emph{space of components $\comp$} of a (topological) space: the set of components endowed
with its natural topology.
We define the \emph{pro-\'etale homotopy type} of a scheme $S$ as the pro-object
\begin{equation*}
    \Pi_{\pe}(S)\colon \Ho(\Hy(S)) \rightarrow \Ho(\sTopo),\; U_{\blt} \mapsto \comp(U_{\blt})
\end{equation*}
in the homotopy category of simplicial spaces.
Here the homotopy categories are defined by identifying (simplicially) homotopic morphisms.
Our definition works for arbitrary schemes, but we have to work with a considerably more complicated object and
lose the intrinsic characterisation of the set of components of a locally connected site.
\par
In \cite{bs14} Bhargav Bhatt and Peter Scholze introduce the notion of w-contractible scheme, which are affine schemes
that trivialise pro-\'etale coverings.
They turn out to be the affine weakly contractible objects in the pro-\'etale site and in Corollary~(\ref{cor:peht-of-w-c})
we show that their pro-\'etale homotopy type is trivial.
In particular, this applies to an algebraically closed field.

\begin{srthm}[Corollary~(\ref{cor:peht-of-w-c})]
Given a w-contractible scheme $W$, there is a natural isomorphism
\begin{equation*}
    \Pi_{\pe}(W) = \const(\comp(W))_{\blt},
\end{equation*}
where $\const$ denotes the constant simplicial object.
\end{srthm}

Moreover, we prove an adjusted version of the statement hinted at in \cite[Remark 4.2.9]{bs14}: the pro-\'etale homotopy
type of a qc scheme can be computed from a single split affine weakly contractible hypercovering.

\begin{srthm}[Theorem~(\ref{thm:sawc-hcov-ini})]
For any split affine weakly contractible hypercovering $W_{\blt}$ in the pro-\'etale site of a qc scheme~$S$ there is a
natural isomorphism
\begin{equation*}
    \Pi_{\pe}(S) = \comp(W_{\blt}).
\end{equation*}
\end{srthm}

We give some evidence that our definition deserves the name of pro-\'etale homotopy type.
First, we prove a statement similar to that of Eric M.~Friedlander in \cite[Proposition 5.2]{fri82}, which relates the set
of simplicial components of the \'etale homotopy type to the set of components of the underlying base scheme.

\begin{srthm}[Theorem~(\ref{thm:peht-comp0})]
Given a scheme $S$, the space of simplicial components of the pro-\'etale homotopy type is given by the space of components
of $S$.
More precisely, we have a natural isomorphism
\begin{equation*}
    \comp \circ \comp_0 \circ \Pi_{\pe}(S) = \comp(S),
\end{equation*}
where $\comp_0$ denotes the space of simplicial components of a simplicial space; see Definition~(\ref{def:simp-top-comp}).
\end{srthm}

Second, we define the cohomology of the pro-\'etale homotopy type of a qcqs scheme with values in pro-\'etale sheaves that
come from sheaves on the space of components; this includes all constant sheaves associated to totally disconnected spaces.
Motivated by their defining property we call these sheaves \emph{$\comp$-sheaves}.

\begin{srthm}[Theorem~(\ref{thm:coho-peht-comp-sh-pe-coho})]
For any $\comp$-sheaf $\sF$ on the pro-\'etale site of a qcqs scheme $S$ there is a natural isomorphism
\begin{equation*}
    \Homol^p(\Pi_{\pe}(S), \sF) = \Homol^p(S_{\pe}, \sF),
\end{equation*}
where the left hand side denotes the cohomology of the pro-\'etale homotopy type and the right hand side the cohomology
of the pro-\'etale site.
\end{srthm}

In section~\ref{sec:peht-field} we make some progress towards describing the pro-\'etale homotopy type of arbitrary fields.
Although we ultimately fall short of our goal, we use our approach in Example~(\ref{exa:peht-reals}) to compute the pro-\'etale
homotopy type of the real numbers.
\par
Lastly, we extend the space of components to pro-objects of simplicial sheaves of the pro-\'etale site and use the theory
by Ilan Barnea and Tomer M.~Schlank introduced in \cite{bs16} to define a refined pro-\'etale homotopy type.
The central result is the following.

\begin{srthm}[Theorem~(\ref{thm:comp-quil-adj})]
Given a qcqs scheme $S$, the extension $\comp_{\sh}$ of the space of components to the pro-\'etale topos together with its
right adjoint $\comp^{\sh}$ induces a pair
\begin{equation*}
    \Pro(\comp_{\sh}) \dashv \Pro(\comp^{\sh})\colon \Pro(\simp\Sh(S_{\pe})) \rlto \Pro(\sCond)
\end{equation*}
of Quillen-adjoint functors with respect to the model structure induced by the natural weak fibration categories on the
pro-\'etale topoi; see \cite[Section 7.1]{bs16}.
\end{srthm}

The \emph{refined pro-\'etale homotopy type} is now given by the value of the left derived functor induced by $\Pro(\comp_{
\sh})$ at the constant pro-(simplicial sheaf) associated to $S$.
We prove that this is indeed a refinement of the pro-\'etale homotopy type defined using hypercoverings.

\begin{srthm}[Theorem~(\ref{thm:peht-vs-refd-peht})]
Let $S$ be a qcqs scheme.
There is a natural isomorphism between the refined pro-\'etale homotopy type and the pro-\'etale homotopy type $\Pi_{\pe}(S)$
considered as pro-objects in the homotopy category of simplicial condensed sets.
\end{srthm}

We also give a definition for the $n$th pro-\'etale homotopy group of a pointed qcqs scheme.
The question of how this notion relates to the pro-\'etale fundamental group turned out to be out of the scope of this paper.

\chapter*{Relation to the literature}

The classical construction of the \'etale homotopy type in \cite{am69} as a pro-object in the homotopy category of simplicial
sets is due to Michael Artin and Barry Mazur.
A refined definition as a pro-object of simplicial sets was given by Eric M.~Friedlander in \cite{fri82}.
The most important application of this theory is in the proof of the Adams conjecture by Dennis Sullivan in \cite{sul74}.
\par
Using a more modern language, Daniel C.~Isaksen considers the \'etale homotopy type from the perspective of model structures
on simplicial presheaves in \cite{isa04}.
In \cite{hoy18} and \cite{car21} Marc Hoyois and David Carchedi define a generalised \'etale homotopy type using the shape
of $\infty$-topoi.
\par
The introduction of the pro-\'etale site by Bhargav Bhatt and Peter Scholze in \cite{bs14} raises the question about a
pro-\'etale homotopy type.
In an appendix to \cite{hrs23} Tamir Hemo, Timo Richarz and Jakob Scholbach use the shape of $\infty$-topoi to define a
pro-\'etale homotopy type and discuss some of its properties.
Furthermore, the group of Clark E.~Barwick, Peter J.~Haine, Tim Holzschuh, Marcin Lara, Catrin Mair, Louis~Martini and Sebastian
Wolf is working on a paper on the pro-\'etale homotopy type using the language of $\infty$-categories and Condensed Mathematics.
Our definitions do not need the machinery of $\infty$-categories: we first follow the approach of Michael Artin and Barry
Mazur and then refine our definition using the language of model categories in the last section.

\chapter*{Structure of the paper}

After briefly explaining some preliminaries and notation in section~\ref{sec:prelim-not}, we begin with reviewing the
notions of w-contractible scheme and weakly contractible object in section~\ref{sec:wc-objs}.
In section~\ref{sec:comp-spc} we introduce the space of components and prove that the pro-\'etale site is not locally
connected.
This enables us to define the pro-\'etale homotopy type in section~\ref{sec:def-peht}.
We then turn towards (affine) weakly contractible hypercoverings and show that they compute the pro-\'etale homotopy type
in section~\ref{sec:wc-hcovs}.
More importantly, we see in section~\ref{sec:peht-from-sawc-hcov} that a single split affine weakly contractible hypercovering
determines the pro-\'etale homotopy type of a qc scheme.
In section~\ref{sec:coho-peht} we define the cohomology of the pro-\'etale homotopy type for a certain class of pro-\'etale
sheaves.
Our attempt to describe the pro-\'etale homotopy type of a field and the computation for the real numbers can be found in
section~\ref{sec:peht-field}.
Lastly, we define a refined version of the pro-\'etale homotopy type in section~\ref{sec:refd-peht} and compare it to our
previous definition.

\chapter*{Acknowledgements}

First, I want to thank Marcin Lara for the helpful and lengthy discussions and his patient answers to my questions.
I want to thank Timo Richarz, Jakob Scholbach and Peter Scholze for answering my questions concerning their respective papers.
Furthermore, I am grateful to Johan Commelin, Bj{\o}rn Dundas and Sebastian Goette for taking their time to answer my questions.
\par
Lastly, but most importantly, I want to thank Annette Huber-Klawitter for her support and help throughout my work to complete
this paper.
It cannot be overstated that this project would not have been possible without her.

\chapter{Preliminaries and notation}\label{sec:prelim-not}

In this section we mention some preliminary concepts and clarify some notation.
We denote the categories of (topological) spaces by $\Topo$, profinite sets by $\Prof$ and extremally disconnected sets
by $\Extr$.
We do not pose any general assumptions on the schemes we consider.
We use the abbreviations \emph{qc} and \emph{qs} for quasi-compact and quasi-separated respectively.
A \emph{pointed scheme} is a scheme together with a geometric point.
\par
For an introduction to sites and sheaves we recommend \cite{tam94} and follow the notation therein.
In particular, a morphism $f\colon \sC \rto \sD$ of sites is given by a functor in the same direction that preserves coverings
and fibre products.
It induces an pair $f_{\sh} \dashv f^{\sh}\colon \Sh(\sC) \rlto \Sh(\sD)$ of adjoint functors between the categories of
sheaves.
We define a \emph{point of a site $\sC$} that admits finite limits to be a morphism $\sC \rto \Set$ of sites that preserves
finite limits and coproducts.
\par
A \emph{weakly \'etale} morphism of schemes is a flat morphism whose diagonal is flat.
For any scheme $S$ its \emph{pro-\'etale site $S_{\pe}$} is defined as the category of weakly \'etale schemes over $S$ together
with families
$\{U_i \rto X \}_{i\in I}$ of morphisms such that any open affine in $X$ is mapped onto by an open affine in $\coprod_{i\in
I} U_i$ as coverings.
We use the abbreviation \emph{covering} for coverings in the pro-\'etale site of a scheme.
We refer the reader to the original source \cite[Section 4]{bs14} or alternatively \cite[Tag 0965]{sta} for a detailed treatment.
\par
In general, we use the definitions and notation of \cite{gj09} when it comes to simplicial objects in a category $\sC$ and
homotopy.
Following \cite{sta} we call the natural functor $\sk_n\colon \simp\sC \rto \simp_n\sC$ from the category of simplicial
objects to their $n$-truncated counter parts the \emph{$n$-skeleton} and its left adjoint $\cosk_n$ the \emph{$n$-coskeleton}.
By abuse of notation we denote the composite $\cosk_n \circ \sk_n$ again by $\cosk_n$.
Apart from section~\ref{sec:refd-peht} the \emph{homotopy category} $\Ho(\simp\sC)$ of simplicial objects in~$\sC$ is the
category with the same objects and equivalence classes of homotopic morphisms.
\par
A \emph{pointed simplicial space} is a simplicial space $X_{\blt}$ and a morphism $\Delta[0] \rto X_{\blt}$.
Together with morphisms that respect the base points they form a category $\sTopo_*$ and their homotopy category
is defined using non-pointed homotopies.
\par
Lastly, a \emph{pro-object in a category $\sC$} is a functor from a variable cofiltered category into $\sC$.
See \cite[Section A.2]{am69} for an introduction including a definition of the morphisms between them.
Recall that a category is cofiltered if for any two objects there is a third object mapping to them and for any two parallel
morphisms there is a third morphism that equalises them; see~\cite[211]{mac98}.
The most important property of pro-objects is that they do not change (up to natural isomorphism) when precomposed with an
initial functor; see \cite[Corollary A.2.5]{am69}.
The central example of an initial functor is a functor $f\colon J \rto I$ between cofiltered categories such that for each
object~$i$ in $I$ there exists an object $j$ in $J$ and a morphism $f(j) \rto i$.
A precise definition of the dual notion can be found in \cite[Paragraph A.1.5]{am69}.
\par
In this paper we encounter set-theoretic issues, which, in the literature, are usually fixed by choosing a universe or an
uncountable strong limit cardinal.
We ignore such issues.

\chapter{Weakly contractible objects}\label{sec:wc-objs}

\begin{srsum}
In this section we review the notions of w-contractible scheme and weakly contractible object.
As an example we describe the affine weakly contractible objects in the pro-\'etale site of an algebraically closed field.
\end{srsum}

\begin{srdef}[{\cite[Definition 1.4]{bs14}}]
An affine scheme $W$ is \emph{w-contractible} if any surjective weakly \'etale morphism $A \rto W$ admits a section.
\end{srdef}

Any affine scheme is covered by a w-contractible one.

\begin{srlem}[{\cite[Lemma 2.4.9]{bs14}}]\label{lem:w-c-covs}
Given an affine scheme $A$, there is a surjective pro-\'etale morphism $W \rto A$ from a w-contractible affine scheme.
\end{srlem}

There is a related relative notion for objects in the pro-\'etale site of a scheme.

\begin{srdef}[{\cite[Tag 090L]{sta}}]
An object $W$ in the pro-\'etale site of a scheme is called \emph{weakly contractible} (or \emph{wc}) if evaluation at
$W$ sends surjections of sheaves to surjections of sets.
\end{srdef}

\begin{srprp}\label{prp:char-awc}
An affine object $W$ in the pro-\'etale site of a scheme is wc if and only if the underlying scheme is w-contractible.
\end{srprp}

\begin{proof}
Can be proved in the same way as \cite[Tag 098H]{sta}.
\end{proof}

We show that affine wc objects behave well with respect to sheafification.

\begin{srlem}\label{lem:shfctn-awc-objs}
For any presheaf $\sF$ and affine wc object $W$ in the pro-\'etale site of a scheme there is a natural isomorphism $\#(\sF)(W)
= \sF(W)$.
Here the symbol $\#$ denotes sheafification.
\end{srlem}

\begin{proof}
By \cite[46]{tam94} the sheafification of $\sF$ is given by applying the $0$th \v{C}ech cohomology of~$W$ twice.
This in turn is defined as a filtered colimit indexed by all coverings $\{U_i \rto W\}_{i\in I}$ over the equaliser of
\begin{equation*}
    \prod_{i\in I} \sF(U_i) \rrto \prod_{i,j\in I} \sF(U_i \times_W U_j).
\end{equation*}
Since $W$ is affine, any covering of $W$ in the pro-\'etale site induces a covering by a single affine object by definition.
We note that the colimit does not change if we only consider coverings of this type.
Moreover, any covering of $W$ by an affine object admits a section since $W$ is affine wc and we conclude that the $0$th
\v{C}ech cohomology is equal to $\sF(W)$.
\end{proof}

An affine wc object in the pro-\'etale site of a scheme is a kind of projective object with respect to coverings.

\begin{srcor}\label{cor:char-awc-proj}
An affine object $W$ in the pro-\'etale site of a scheme $S$ is wc if and only if for any diagram
\begin{equation*}
\begin{tikzpicture}
    \node (W) at (0, 0)   {$W$};
    \node (U) at (3, 1.5) {$U$};
    \node (X) at (3, 0)   {$X\hsmash{,}$};
    \path[-{Straight Barb[length=2pt]}, font=\scriptsize, line width=0.2mm]
        (W) edge (X)
        (U) edge (X);
    \path[dashed, -{Straight Barb[length=2pt]}, font=\scriptsize, line width=0.2mm]
        (W) edge (U);
\end{tikzpicture}
\end{equation*}
where $U \rightarrow X$ is a covering, there is a dashed morphism making the diagram commute.
\end{srcor}

\begin{proof}
Assume that $W$ is an affine wc object in the pro-\'etale site of a scheme and that we are given a diagram as in the statement.
The base change $W \times_X U \rto W$ is a covering and since $W$ is qc, it is refined by an affine covering $A \rto W$.
By the definition of affine wc objects we obtain a section $W \rto A$ and thus a morphism $W \rto W \times_X U \rto U$.
This is the desired lift.
\par
Assume that $W$ is an affine object in the pro-\'etale site of a scheme that satisfies the lifting property above.
Given an affine covering $B \rto W$ we may apply the lifting property to this covering and the identity on $W$ to obtain a
section $W \rto B$.
\end{proof}

Using this characterisation we can determine the affine wc objects in the pro-\'etale site of an algebraically closed
field.
Recall the notion of extremally disconnected space.

\begin{srdef}[{\cite[Tag 08YI]{sta}}]\label{def:ex-disc}
A space is called \emph{extremally disconnected} if it is profinite and the closure of every open subset is open.
\end{srdef}

The following result is stated in \cite[Remark 4.1.11]{bs14}.

\begin{srprp}\label{prp:awc-pe-site-k-alg}
Under the equivalence from \cite[Example 4.1.10]{bs14} the affine wc objects in the pro-\'etale site of an algebraically
closed field correspond to extremally disconnected spaces.
\end{srprp}

\begin{proof}
Using Corollary~(\ref{cor:char-awc-proj}) we see that the affine wc objects in the pro-\'etale site of $k$ correspond to
the projective objects in the category of profinite sets.
By \cite[Theorem 2.5]{gle58} the projective objects in the category of profinite sets are exactly the extremally
disconnected ones.
\end{proof}

Lastly, we prove that any (pointed) covering of an object in the pro-\'etale site of a scheme can be refined by a (pointed)
wc covering.

\begin{srlem}\label{lem:cov-refd-by-wc}
Any covering $U \rto X$ in the pro-\'etale site of a scheme has a refinement by a wc covering $W \rto X$ and if $X$ is qc,
then $W$ may be chosen affine.
Moreover, the same is true if all the involved schemes are pointed.
\end{srlem}

\begin{proof}
By the definition of coverings in the pro-\'etale site we obtain a covering $\{U_i \rto X\}_{i\in I}$ by open affine subschemes
of $U$.
In the pointed case we may assume that there is a $U_{i_0}$ that contains the point $u$.
For each $U_i$ we choose affine wc coverings $W_i$ and consider their coproduct $W$.
By construction the natural morphism $W \rto X$ is a covering and it is easy to see that $W$ is wc as a coproduct of wc
objects.
In the pointed case we may choose a preimage of the point $u$ in $W_{i_0}$.
If the object $X$ is qc there is a finite subcovering of the covering $\{U_i \rto X\}_{i\in I}$.
\end{proof}

\chapter{Space of components}\label{sec:comp-spc}

\begin{srsum}
We recall the structure of site on the space of components of a qcqs scheme and the related morphism into the pro-\'etale
site of the underlying scheme.
Using this result we show that any qcqs scheme admits pro-\'etale coverings from schemes that are not locally connected.
\end{srsum}

\begin{srdef}
For any point $x$ of a space $X$, the \emph{component of $x$} is the maximal, with respect to inclusion, connected
subspace of $X$ that contains $x$.
The \emph{space of components of $X$} is the set of components, denoted by $\comp(X)$, endowed with the quotient
topology induced by the natural surjection $X \rightarrow \comp(X)$.
\end{srdef}

Moreover, recall the notion of totally disconnected space.

\begin{srdef}[{\cite[Tag 04MC]{sta}}]
A space is called \emph{totally disconnected} if all its components consist of a single point.
\end{srdef}

\begin{srrem}
The space of components is totally disconnected by \cite[Tag 08ZL]{sta} and gives rise to a functor $\comp\colon \Topo
\rightarrow \TotDisc$.
\end{srrem}

The space of components determines morphisms of the underlying space into totally disconnected spaces.

\begin{srprp}[{\cite[Tag 08ZL]{sta}}]\label{prp:comp-incl-adj}
The space of components and the inclusion $\iota$ of the full subcategory of totally disconnected spaces form a pair
\begin{equation*}
    \comp \dashv \iota\colon \Topo \rightleftarrows \TotDisc.
\end{equation*}
of adjoint functors.
\end{srprp}

The space of components of a qcqs scheme is particularly nice.
We learnt about the following argument from Marcin Lara.

\begin{srlem}[{\cite[30]{bs14}}]\label{lem:comp-of-qcqs-is-prof}
For a qcqs scheme $X$ its space of components $\comp(X)$ is profinite.
\end{srlem}

\begin{proof}
By \cite[Corollaire 8.5]{laz67} the component of a point of $X$ equals the intersection of all closed-open neighbourhoods
of the point.
We conclude the proof by applying \cite[Tag 0900]{sta} since $X$ is qc.
\end{proof}

In the \'etale setting the question whether the set of components should be considered with its natural topology does
not pose itself.

\begin{srrem}\label{rem:etale-over-loc-noeth}
The components of a locally connected space (e.g.\ a locally Noetherian scheme) are open by
\cite[Proposition I.11.6.11]{bou07} and thus its space of components is discrete.
The property of being locally Noetherian lifts along \'etale morphisms by \cite[Tag 01T6]{sta}.
Since the base scheme in the setting of \cite{am69} and \cite{fri82} is locally Noetherian, all schemes are locally connected.
The same is not true in the pro-\'etale setting.
\end{srrem}

In the following we use the space of components to construct objects in the pro-\'etale site that are not locally
connected.

\begin{srdef}[{\cite[30]{bs14}}]
For any qcqs scheme $S$ the \emph{pro-\'etale site $\comp(S)_{\pe}$ on the space of components} is defined to be the
category $\Prof/\comp(S)$ together with families of morphisms that have a finite surjective subfamily as coverings.
\end{srdef}


Using the following lemma the space of components induces a natural morphism of sites.

\begin{srlem}[{\cite[10]{bs14}}]\label{lem:comp-fibr-prod-adj}
Let $A$ be an affine scheme.
For any profinite set $P$ over $\comp(A)$ the fibre product $\alpha_A(P) \coloneqq P \times_{\comp(A)} A$ of spaces
admits a natural structure of an affine scheme.
Together with the space of components on the subcategory $A_{\aff,\pe}$ of affine objects in the pro-\'etale site of $A$
it induces a pair
\begin{equation*}
    \comp \dashv \alpha_A\colon A_{\aff,\pe} \rightleftarrows \Prof/\comp(A)
\end{equation*}
of adjoint functors.
\end{srlem}

As stated in \cite[30]{bs14} we extend this construction to the pro-\'etale site of a qcqs scheme.

\begin{srprp}[{\cite[30]{bs14}}]\label{prp:con-pe-covs}
For any qcqs scheme $S$ and profinite set $P$ over $\comp(S)$ the fibre product $P \times_{\comp(S)} S$ of spaces admits
a natural scheme-structure by pulling back the structure sheaf of $S$.
This induces a morphism $\alpha\colon \comp(S)_{\pe} \rightarrow S_{\pe}$ of sites.
\end{srprp}

\begin{proof}
Let $P$ be a profinite set over $\comp(S)$.
We endow the fibre product $P \times_{\comp(S)} S$ with the structure of a locally ringed space by pulling back the
structure sheaf of $S$ along the projection to $S$.
Given an open affine subscheme $U$ of $S$, we obtain an open subspace
\begin{equation*}
    (P \times_{\comp(S)} \comp(U)) \times_{\comp(U)} U = P \times_{\comp(S)} U
\end{equation*}
of the space $P \times_{\comp(S)} S$, which is an affine scheme by Lemma~(\ref{lem:comp-fibr-prod-adj}).
Hence any open affine cover of $S$ gives rise to an open affine cover of $P \times_{\comp(S)} S$ and we conclude that $P
\times_{\comp(S)} S$ is a scheme and weakly \'etale over $S$.
\par
We deduce the functoriality in the profinite set $P$ of the constructed scheme from the functoriality of the adjunction
of Lemma~(\ref{lem:comp-fibr-prod-adj}).
The functor $\alpha$ preserves coproducts, singleton coverings and fibre products and is therefore a morphism of sites.
\end{proof}

The following lemma is proved in \cite[Lemma 2.2.8]{bs14} and \cite[Tag 096C]{sta} in the affine case.
We use the same proof to extend the result to qcqs schemes.

\begin{srlem}\label{lem:comp-pe-covs}
For any qcqs scheme $S$ and profinite set $P$ there is a natural isomorphism
\begin{equation*}
    \comp(P \times_{\comp(S)} S) = P,
\end{equation*}
where the fibre product is in the category of spaces.
\end{srlem}

\begin{proof}
Since the space of components preserves surjections the morphism $\comp(P \times_{\comp(S)} S) \rightarrow P$
is surjective.
Let $p$ be a point in $P$ and $s$ its image in $\comp(S)$.
We see that the fibre of $p$ under the projection $P \times_{\comp(S)} S \rightarrow P$ is isomorphic to the fibre of $s$
under the morphism $S \rightarrow \comp(S)$, which is connected.
This implies that the morphism $\comp(P \times_{\comp(S)} S) \rightarrow P$ is a bijection and therefore an isomorphism
by \cite[Tag 08YE]{sta}.
\end{proof}

Using this construction we can construct new w-contractible schemes from a given one.

\begin{srlem}\label{lem:con-w-c}
Given a w-contractible scheme $W$ and a morphism $E \rto \comp(W)$ from an extremally disconnected space, the scheme $E
\times_{\comp(W)} W$ is w-contractible.
\end{srlem}

\begin{proof}
We use the characterisation of \cite[Theorem 1.8]{bs14}.
The space $\alpha(W)\coloneqq E \times_{\comp(W)} W$ is w-local by \cite[Tag 096C]{sta}.
Given a closed point $p$ of $\alpha(W)$ its image $w$ is closed in $W$ and we compute $\sO_{\alpha(W), p} = \sO_{W, w}$,
which is strictly Henselian.
Moreover, the space of components of $\alpha(W)$ is the extremally disconnected space $E$ by Lemma~(\ref{lem:comp-pe-covs}).
\end{proof}

We show that any qcqs scheme admits a pro-\'etale covering from a scheme that is not locally connected.

\begin{srexa}\label{exa:we-cov-from-non-lc}
Let $S$ be a qcqs scheme and $E \rightarrow \comp(S)$ a surjection from an infinite discrete set.
Using the Stone-\v{C}ech compactification and its universal property we obtain a surjective morphism $\upbeta E \rightarrow
\comp(S)$ of spaces from a non-discrete profinite set.
In fact, the space $\upbeta E$ is qc and contains $E$ as an open subspace by~\cite[Tag 090A]{sta}; therefore it can not
be discrete.
By the universal property of the Stone-\v{C}ech compactification $\upbeta E$ is a projective object in the category of
qc Hausdorff spaces and therefore extremally disconnected by \cite[Theorem~2.5]{gle58}.
We obtain a weakly \'etale covering $\upbeta E \times_{\comp(S)} S \rightarrow S$ using Proposition~(\ref{prp:con-pe-covs}),
which is not locally connected.
\end{srexa}

For an arbitrary scheme we still obtain a pro-\'etale morphism from a scheme that is not locally connected.

\begin{srexa}\label{exa:pe-morph-from-non-lc}
Let $U$ be an open affine subscheme of a scheme $S$.
We apply the construction of Example~(\ref{exa:we-cov-from-non-lc}) to the scheme $U$ and obtain a pro-\'etale morphism
$\upbeta E \times_{\comp(U)} U \rightarrow U$ by Lemma~(\ref{lem:comp-fibr-prod-adj}).
By composing with the open immersion $U \rightarrow S$ we get a pro-\'etale morphism to~$S$ from a scheme that is not locally
connected.
\end{srexa}

Lastly, we introduce some notation for pointed spaces.

\begin{srrem}
Given a pointed space $(X,x)$, the space of components $\comp(X)$ is naturally pointed in the component that contains
the point $x$.
By abuse of notation we will denote this component again by the symbol $x$.
We obtain a functor $\comp\colon \Topo_* \rto \TotDisc_*$.
\end{srrem}

\chapter{Definition of the pro-\'etale homotopy type}\label{sec:def-peht}

\begin{srsum}
In this section we define the pro-\'etale homotopy type of a (pointed) scheme and prove a analogue of
\cite[Proposition 5.2]{fri82}.
\end{srsum}

We recall the notion of pointed simplicial space and pointed hypercovering.

\begin{srrem}
A hypercovering of a scheme $S$ is a simplicial object $U_{\blt}$ of its pro-\'etale site such that the natural morphisms
$U_{n+1} \rto \cosk_n(U_{\blt})_{n+1}$ as well as $U_0 \rto S$ are coverings.
Together with morphisms of simplicial objects they form a category $\Hy(S)$ and their homotopy category obtained by identifying
homotopic morphisms is cofiltered by \cite[Th\'eor\`eme V.7.3.2.1]{gro72b}.
\par
A geometric point $s$ of $S$ induces a point $s$ of the pro-\'etale site of $S$ by mapping schemes to their geometric points
over $s$.
A \emph{pointed hypercovering} is a hypercovering $U_{\blt}$ of  $S$ and a morphism $\Delta[0] \rto s(U_{\blt})$.
Together with morphisms respecting the base points we obtain a category $\Hy(S, s)$ of pointed hypercoverings.
As above their homotopy category is defined by identifying homotopic morphisms and cofiltered.
Note that we use non-pointed homotopies since this does not affect the definition of the pro-\'etale homotopy type.
\end{srrem}

We are not able to use the definition of the homotopy type introduced in \cite[112]{am69} since the pro-\'etale site is
not locally connected.

\begin{srrem}\label{rem:pe-site-not-loc-con}
As Example~(\ref{exa:pe-morph-from-non-lc}) shows, any scheme admits a weakly \'etale morphism from a scheme that is not
locally connected.
This is a consequence of weakly \'etale morphisms not possessing the finiteness condition of \'etale morphisms.
Therefore we take the natural topology on the set of components into account.
\end{srrem}

We introduce the central notion in its pointed and non-pointed flavour.

\begin{srdef}\label{def:peht}
For any pointed scheme $(S,s)$ its \emph{(pointed) pro-\'etale homotopy type} is the pro-object
\begin{equation*}
    \Pi_{\pe}(S,s)\colon \Ho(\Hy(S,s)) \rightarrow \Ho(\sTopo_*),\; (U_{\blt},u) \mapsto \comp(U_{\blt},u)
\end{equation*}
in the homotopy category of pointed simplicial spaces.
An analogous definition using non-pointed hypercoverings gives rise to the \emph{pro-\'etale homotopy type $\Pi_{\pe}(S)$ of $S$}.
\end{srdef}

We prove that both definitions are compatible.

\begin{srlem}
For a pointed scheme $(S, s)$, the natural functor $\Ho(\Hy(S, s)) \rightarrow \Ho(\Hy(S))$ is initial.
Moreover, there is a natural isomorphism
\begin{equation*}
    \Pi_{\pe}(S) = \Pi_{\pe}(S,s)
\end{equation*}
of pro-objects in $\Ho(\sTopo))$.
\end{srlem}

\begin{proof}
We conclude the first claim from the fact that geometric points lift along coverings.
The second claim follows from \cite[Corollary A.2.5]{am69}.
\end{proof}

As a first test, we want to prove a version of \cite[Proposition 5.2]{fri82}: composing the \'etale homotopy type with the
set of simplicial components gives back the set of components of the base scheme.
First, we introduce an analogue of the set of components for simplicial spaces.

\begin{srdef}\label{def:simp-top-comp}
Given a simplicial space $X_{\blt}$, its \emph{space of (simplicial) components} is the space
\begin{equation*}
    \comp_0(X_{\blt}) \coloneqq \CoEq(X_1 \rightrightarrows X_0).
\end{equation*}
Denote the induced functor by $\comp_0\colon \sTopo \rightarrow \Topo$.
\end{srdef}

We can characterise the component functor of simplicial spaces as an adjoint.

\begin{srlem}\label{lem:comp-const-adj}
The component functor $\comp_0$ of simplicial spaces and the constant simplicial space functor $\const$ form a pair
\begin{equation*}
    \comp_0 \dashv \const\colon \sTopo \rightleftarrows \Topo
\end{equation*}
of adjoint functors.
\end{srlem}

\begin{proof}
Let $\comp_0(X_{\blt}) \rightarrow Y$ be a morphism of spaces.
We construct a morphism $X_{\blt} \rightarrow \const(Y)_{\blt}$ of simplicial spaces inductively.
Composition with the natural morphism $X_0 \rightarrow \comp_0(X_{\blt})$ gives rise to a morphism $X_0 \rightarrow Y$,
which coequalises the face morphisms $d_0$ and $d_1$.
\par
Suppose a morphism $f_n\colon X_n \rightarrow Y$ that coequalises the $d_i$'s of degree $n{+}1$ has been constructed.
This induces a morphism $f_{n+1}\colon X_{n+1} \rightarrow Y$ by pre-composing with any of these $d_i$'s.
Using the induction hypothesis and the simplicial identities one computes that $f_{n+1}$ coequalises the face morphisms
of degree $n{+}2$.
\par
A morphism $X_{\blt} \rightarrow \const(Y)_{\blt}$ of simplicial spaces induces a unique morphism $\comp_0(X_{\blt})
\rightarrow Y$ by the universal property of the coequaliser.
These constructions are mutually inverse.
\end{proof}

Analogous to the result in \cite[Tag 00GX]{ker} for simplicial sets, we prove that the space of components of a
simplicial space preserves products.

\begin{srlem}\label{lem:comp-presrv-prod}
Given simplicial spaces $X_{\blt}$ and $Y_{\blt}$, there is a natural isomorphism
\begin{equation*}
    \comp_0(X_{\blt} \times Y_{\blt}) = \comp_0(X_{\blt}) \times \comp_0(Y_{\blt}).
\end{equation*}
\end{srlem}

\begin{proof}
Applying the universal property of the coequaliser to $X_0 \times Y_0 \rightarrow X_0 \rightarrow \comp_0(X_{\blt})$ yields
a morphism of spaces $\comp_0(X_{\blt} \times Y_{\blt}) \rto \comp_0(X_{\blt})$.
The same works for $\comp_0(Y_{\blt})$ and we obtain a morphism of spaces into their product, which is a bijection by
\cite[Tag 00GX]{ker}.
Since the morphism $X_0 \times Y_0 \rto \comp_0(X_{\blt} \times Y_{\blt})$ is surjective and the morphism $X_0 \times Y_0
\rto \comp_0(X_{\blt}) \times \comp_0(Y_{\blt})$ is open, the bijective morphism constructed above is open by
\cite[Proposition I.5.1.1.b]{bou07}.
\end{proof}

We are now able to conclude that the component functor of simplicial spaces identifies homotopic morphisms.

\begin{srcor}\label{cor:scomp-ho-stop}
The space of components induces a functor~$\Ho(\sTopo) \rightarrow \Topo$, which is again denoted by $\comp_0$.
\end{srcor}

\begin{proof}
Given a simplicial space $X_{\blt}$, we compute
\begin{equation*}
    \comp_0(X_{\blt} \times \Delta[1]) = \comp_0(X_{\blt}) \times \comp_0(\Delta[1])
        = \comp_0(X_{\blt}).
\end{equation*}
Since the natural isomorphism on the right is induced by projection, composing with the natural morphisms $\comp_0(X_{\blt})
\rto \comp_0(X_{\blt} \times \Delta[1])$, induced by the inclusions $\iota_j\colon \Delta[0] \rightarrow \Delta[1]$ for
$j = 1,2$, gives back the identity.
Therefore any two morphisms of simplicial spaces related by a homotopy become equal after applying $\comp_0$.
\end{proof}

We introduce some notation for pointed simplicial spaces.

\begin{srrem}
The space of components of a pointed simplicial space $(X_{\blt}, x)$ is naturally pointed and hence we obtain a functor
$\comp_0\colon \Ho(\sTopo_*) \rightarrow \Topo_*$ by Corollary~(\ref{cor:scomp-ho-stop}).
\end{srrem}

Our goal is to prove that applying this functor to the pro-\'etale homotopy type gives back the space of components of
the underlying scheme.
The original proof makes use of the fact that the functor $\Hom_{\Topo}(\earg, S)$ is a sheaf on the \'etale site.
We use the following generalisation of this; another sign that the pro-\'etale site interacts well with topological objects.

\begin{srlem}[{\cite[Lemma 4.2.12]{bs14}}]\label{lem:sheaf-assoc-to-space}
For any scheme $S$ and space $Y$ the functor
\begin{equation*}
    \ul{Y}\colon S_{\pe} \rightarrow \Set,\; X \mapsto \Hom_{\Topo}(X, Y)
\end{equation*}
defines a sheaf on $S_{\pe}$, called the \emph{constant sheaf associated to $Y$}.
\end{srlem}

Already for a point we obtain interesting objects.

\begin{srexa}
By \cite[Example 4.1.10]{bs14} the sheaves on the pro-\'etale site of an algebraically closed field are equivalent to the
sheaves on the canonical site of profinite sets and also known as condensed sets.
The constant sheaf associated to a space~$Y$ is nothing but the condensed set associated to $Y$ by \cite[Example 1.5]{sch19}.
\end{srexa}

We now generalise \cite[Proposition 5.2]{fri82}.

\begin{srthm}\label{thm:peht-comp0}
For any pointed scheme $(S,s)$ there is a natural isomorphism
\begin{equation*}
    \comp \circ \comp_0 \circ \Pi_{\pe}(S,s) = \comp(S,s),
\end{equation*}
as pro-objects of pointed spaces.
\end{srthm}

\begin{proof}
We show that for any pointed hypercovering $(U_{\blt}, u)$ and totally disconnected space $T$ there is a natural bijection
\begin{equation*}
    \Hom_{\Topo}((\comp \circ \comp_0 \circ \comp)(U_{\blt}), T) = \Hom_{\Topo}(\comp(S), T),
\end{equation*}
induced by the morphism $U_0 \rightarrow S$.
Since the category of totally disconnected spaces is a full subcategory of the category of spaces, this implies the claim.
\par
Starting with the expression on the left, we use Lemma~(\ref{prp:comp-incl-adj}) and Lemma~(\ref{lem:comp-const-adj})
to calculate
\useshortskip
\begin{align*}
    \Hom_{\Topo}((\comp \circ \comp_0 \circ \comp)(U_{\blt}), T) &= \Hom_{\sTopo}(U_{\blt}, \const(T)_{\blt})\\
        &= \Hom_{\Topo}(\comp_0(U_{\blt}), T)\\
        &= \ulsi{T}(\comp_0(U_{\blt})).
\end{align*}
Since the space $\comp_0(U_{\blt})$ is defined as a coequaliser and the sheaf $\ulsi{T}$ sends colimits to limits, the diagram
\begin{equation*}
    \ulsi{T}(\comp_0(U_{\blt})) \rightarrow \ulsi{T}(U_0) \rightrightarrows \ulsi{T}(U_1)
\end{equation*}
is an equaliser.
Since hypercoverings compute global sections we conclude $\ulsi{T}(\comp_0(U_{\blt})) = \ulsi{T}(S)$.
Another application of Lemma~(\ref{prp:comp-incl-adj}) proves the claim.
\end{proof}

\begin{srrem}
When compared to \cite[Proposition 5.2]{fri82}, we need an additional application of the space of components to guarantee
that the space $\comp \circ \comp_0 \circ \Pi_{\pe}(S,s)$ from Theorem~(\ref{thm:peht-comp0}) is totally disconnected.
We thank Marcin~Lara for pointing out that a quotient of a totally disconnected space need not be totally disconnected.
\end{srrem}

\chapter{Weakly contractible hypercoverings}\label{sec:wc-hcovs}

\begin{srsum}
We show that the pro-\'etale homotopy type of a scheme can be computed by wc hypercoverings.
For a qcqs scheme we show that it is determined by affine wc hypercoverings.
\end{srsum}

We recall the notion of split hypercovering and their construction result.

\begin{srrem}\label{rem:con-hcov}
Following \cite[D\'efinition Vbis.5.1.1]{gro72b} we call a hypercovering \emph{split} if in each degree it decomposes into
a disjoint union of its non-degenerate and degenerate simplices.
Given an $n$-truncated split hypercovering $U_{\blt}$ and a morphism $N \rto \cosk_n(U_{\blt})_{n+1}$, there is a unique
$(n{+}1)$-truncated hypercovering $U'_{\blt}$ that has $N$ as its non-degenerate simplices.
Giving a morphism $U'_{\blt} \rto V_{\blt}$ of $(n{+}1)$-truncated hypercoverings is the same as giving a morphism on
$n$-truncations and a morphism $N \rto V_{n+1}$ compatible with the morphisms to the $n$-coskeleta.
One may consult \cite[Proposition Vbis.5.1.3]{gro72b} or \cite[Section 6.2]{del74} for more details.
\end{srrem}

We introduce some language for pointed hypercoverings.

\begin{srdef}
A pointed hypercovering $(U_{\blt},u)$ of a pointed scheme is called \emph{wc} (\emph{qc}, or \emph{affine}) if all the
$U_n$ are wc (qc, or affine) respectively.
We denote the category of pointed wc (or affine wc) hypercoverings of a pointed scheme $(S,s)$ by $\Hywc(S,s)$ (or $\Hyawc(S,s)$)
respectively.
\end{srdef}

Using split hypercoverings we construct wc refinements of pointed hypercoverings.

\begin{srlem}\label{lem:hcov-refd-by-wc}
Any pointed hypercovering $(U_{\blt}, u)$ of a pointed scheme $(S, s)$ admits a surjective morphism $(W_{\blt}, w) \rto
(U_{\blt}, u)$ from a pointed split wc hypercovering.
If the scheme $S$ is qc there is a morphism $(W'_{\blt}, w') \rto (U_{\blt}, u)$ from a pointed split affine wc hypercovering.
\end{srlem}

\begin{proof}
We choose a pointed wc covering $(W_0, w)$ that refines $(U_0, u)$ using Lemma~(\ref{lem:cov-refd-by-wc}).
This is a surjective morphism from a pointed $0$-truncated split wc hypercovering.
\par
Assume that a surjective morphism $(W_{\blt}, w) \rto (\sk_n(U_{\blt})_{\blt}, u)$ from a pointed $n$-truncated split wc
hypercovering has been constructed.
The projection to $\cosk_n(W_{\blt})_{n+1}$ of the fibre product $U_{n+1} \times_{\cosk_n(U_{\blt})_{n+1}} \cosk_n(W_{\blt}
)_{n+1}$ is a covering.
By Lemma~(\ref{lem:cov-refd-by-wc}) it can be refined by a wc covering $N_{n+1} \rto \cosk_n(W_{\blt})_{n+1}$.
Therefore we obtain an extension to a surjective morphism $(W_{\blt},w) \rightarrow (\sk_{n+1}(
U_{\blt})_{\blt}, u)$ from a pointed $(n{+}1)$-truncated split wc hypercovering by Remark~(\ref{rem:con-hcov}).
\par
For the second claim we use the same technique together with the fact that fibre products of qc schemes and hence the
$n$-coskeleton of qc hypercoverings are qc.
Note that the morphism in degree $0$ need no longer be surjective.
\end{proof}

As a consequence, we conclude that the pro-\'etale homotopy type is determined by pointed wc hypercoverings.

\begin{srcor}\label{cor:wc-hcovs-ini}
For any pointed scheme $(S, s)$ the category $\Ho(\Hywc(S,s))$ is cofiltered and its inclusion into $\Ho(\Hy(S,s))$ is
initial.
Moreover, there is a natural isomorphism
\begin{equation*}
    \left.{\Pi_{\pe}(S,s)}\right|_{\Ho(\Hywc)} = \Pi_{\pe}(S,s)
\end{equation*}
of pro-objects in $\Ho(\sTopo_*))$.
If $S$ is qc the same holds for the category $\Ho(\Hyawc(S, s))$.
\end{srcor}

\begin{proof}
Since the homotopy category of pointed hypercoverings is cofiltered, the homotopy category of pointed wc hypercoverings
is cofiltered and its inclusion is initial by Lemma~(\ref{lem:hcov-refd-by-wc}).
The natural isomorphisms exists by \cite[Corollary A.2.5]{am69} and the same arguments apply in the affine case.
\end{proof}

For a qc scheme we are only dealing with profinite sets.

\begin{srcor}\label{cor:peht-qcqs-sprof}
The pro-\'etale homotopy type of a qc scheme is a pro-object in the homotopy category of pointed simplicial profinite
sets.
\end{srcor}

\begin{proof}
Combine Corollary~(\ref{cor:wc-hcovs-ini}) with Lemma~(\ref{lem:comp-of-qcqs-is-prof}).
\end{proof}

\chapter{Pro-\'etale homotopy type from a single hypercovering}\label{sec:peht-from-sawc-hcov}

\begin{srsum}
We prove that the pro-\'etale homotopy type of a qc scheme is determined by a single split affine wc hypercovering; a
statement hinted at in \cite[Remark 4.2.9]{bs14}.
Additionally, we prove that w-contractible schemes have a trivial pro-\'etale homotopy type.
\end{srsum}

We first state the central result of this section.

\begin{srthm}\label{thm:hcov-sawc-ref}
Given a split affine wc hypercovering $W_{\blt}$ and a hypercovering $U_{\blt}$ in the pro-\'etale site of a qc scheme,
there is a unique morphism $W_{\blt} \rightarrow U_{\blt}$ up to homotopy equivalence.
\end{srthm}

As a direct consequence we obtain that the pro-\'etale homotopy type of a qc scheme can be computed from a single split
affine wc hypercovering.

\begin{srthm}\label{thm:sawc-hcov-ini}
The inclusion of any split affine wc hypercovering $W_{\blt}$ into the homotopy category of all hypercoverings in the pro-\'etale
site of a qc scheme $S$ is initial.
Moreover, there is a natural isomorphism
\begin{equation*}
    \Pi_{\pe}(S) = \comp(W_{\blt}).
\end{equation*}
\end{srthm}

\begin{srrem}
Note that requiring the affine wc hypercovering to be split is not a restrictive assumption.
Let $W_{\blt}$ be an affine wc hypercovering.
By Remark~(\ref{rem:con-hcov}) we obtain a natural split hypercovering that has $W_n$ as its non-degenerate simplices in
degree $n$.
Moreover, it is affine wc since this property is preserved by finite disjoint unions.
\end{srrem}

\begin{srcor}\label{cor:peht-of-w-c}
Given a w-contractible scheme $W$, the inclusion of $\const(W)_{\blt}$ into the homotopy category of all hypercoverings
in the pro-\'etale site of $W$ is initial.
Moreover, there is a natural isomorphism
\begin{equation*}
    \Pi_{\pe}(W) = \const(\comp(W))_{\blt}.
\end{equation*}
\end{srcor}

\begin{proof}
Use the fact that the constant simplicial object $\const(W)_{\blt}$ is a split affine wc hypercovering.
Note that we only need the existence part of Theorem~(\ref{thm:hcov-sawc-ref}) since the only endomorphism of the constant
simplicial object $\const(W)_{\blt}$ is the identity.
\end{proof}

We start by proving the existence part of Theorem~(\ref{thm:hcov-sawc-ref}).

\begin{srlem}\label{lem:hcov-sawc-ref}
Given a split affine wc hypercovering $W_{\blt}$ and a hypercovering $U_{\blt}$ in the pro-\'etale site of a qc scheme,
there is a morphism $W_{\blt} \rightarrow U_{\blt}$.
\end{srlem}

\begin{proof}
We construct a morphism $W_{\blt} \rightarrow U_{\blt}$ inductively using the lifting property of affine wc objects in
the pro-\'etale site.
Since $W_0$ is affine wc, we obtain a morphism $W_0 \rto U_0$ of $0$-truncated split hypercoverings by lifting $W_0 \rto
S$ along the covering $U_0 \rto S$; see Corollary~(\ref{cor:char-awc-proj}).
\par
Assume that a morphism $\sk_n(W_{\blt})_{\blt} \rightarrow \sk_n(U_{\blt})_{\blt}$ of $n$-truncated hypercoverings has
been constructed.
Consider the diagram
\begin{equation*}
\begin{tikzpicture}
    \node (NWn+1) at (0,   0) {$\NDeg W_{n+1}$};
    \node (cWn+1) at (3,   0) {$\cosk_n(W_{\blt})_{n+1}$};
    \node (cUn+1) at (6,   0) {$\cosk_n(U_{\blt})_{n+1}\hsmash{.}$};
    \node (Un+1)  at (6, 1.5) {$U_{n+1}$};
    \path[-{Straight Barb[length=2pt]}, font=\scriptsize, line width=0.2mm]
        (NWn+1) edge (cWn+1)
        (cWn+1) edge (cUn+1)
        (Un+1)  edge (cUn+1);
    \path[dashed, -{Straight Barb[length=2pt]}, font=\scriptsize, line width=0.2mm]
        (NWn+1) edge (Un+1);
\end{tikzpicture}
\end{equation*}
Since the affine wc hypercovering $W_{\blt}$ is split, the summand $\NDeg W_{n+1}$ is affine wc and the dashed lift
exists.
We obtain a morphism $\sk_{n+1}(W_{\blt})_{\blt} \rto \sk_{n+1}(U_{\blt})_{\blt}$ of $(n{+}1)$-truncated hypercoverings
by Remark~(\ref{rem:con-hcov}).
\end{proof}

\begin{srrem}
In the pointed setting we run into issues.
If $(W_{\blt}, w)$ is a pointed split affine wc hypercovering and $(U_{\blt}, u)$ a pointed hypercovering, we obtain a
morphism $W_0 \rightarrow U_0$.
We need to show that this morphism can be chosen to map $w$ to the point $u$.
More specifically, we need to prove that for any two geometric points $w_1$ and $w_2$ of an affine wc object $W$ there is
a morphism~$(W, w_1) \rto (W, w_2)$.
\par
An alternative solution, suggested by Annette Huber-Klawitter, would be to define pointed w-contractible schemes using a
pointed lifting property.
Then the difficulty lies in proving the existence of such pointed w-contractible schemes.
\end{srrem}

Lastly, we have to prove that any endomorphism of a split affine wc hypercovering in the pro-\'etale site of a qc scheme
is homotopic to the identity.

\begin{srrem}
The classical definition of homotopy as in \cite[Definitions 5.1]{may67} does not interact well with coskeleta.
We thank Jakob Scholbach for pointing us towards the slightly adjusted version in \cite[Section 3]{bkr19}, which interacts
well with coskeleta.
This notion of reduced homotopy already appears in \cite[Tag 019L]{sta}.
As far as we are aware the fact that reduced homotopies interact well with coskeleta does not appear in \cite{sta}.
\end{srrem}

\begin{srdef}[{\cite[3]{bkr19}}]
Let $f,g\colon X_{\blt} \rto Y_{\blt}$ be morphisms of simplicial objects in a category.
A \emph{reduced homotopy} from $f$ to $g$ is a family $\{r^i_n\colon X_n \rto Y_n\}_{0 \leqslant i \leqslant n+1}$ of
morphisms for each $n \geqslant 0$ such that $r^0_n = f_n$ and $r^{n+1}_n = g_n$ and
\begin{alignat*}{1}
    d^j_n \circ r^i_n &=
        \begin{cases*}
            r^{i-1}_{n-1} \circ d^j_n\quad &if $i > j$\\
            r^{i}_{n-1} \circ d^j_n\quad &if $i \leqslant j$,
        \end{cases*}\\
    s^j_n \circ r^i_n &=
        \begin{cases*}
            \mathrlap{r^{i+1}_{n+1} \circ s^j_n}\hphantom{r^{i-1}_{n-1} \circ d^j_n}\quad &if $i > j$\\
            r^{i}_{n+1} \circ s^j_n\quad &if $i \leqslant j$.
        \end{cases*}
\end{alignat*}
A \emph{$d$-truncated reduced homotopy} is a family $\{r^i_n\colon X_n \rto Y_n\}_{0 \leqslant i \leqslant n+1}$ for all
$n \leqslant d$ such that the identities are satisfied up to degree $d$.
\end{srdef}

\begin{srprp}[{\cite[Proposition 5.2]{bkr19}}]\label{prp:n-trunc-red-ho-cosk}
Let $f,g\colon X_{\blt} \rto Y_{\blt}$ be morphisms of $n$-truncated simplicial objects in a category that has finite
limits.
Any $n$-truncated reduced homotopy $r^i_{\blt}$ between $f$ and $g$ induces an $(n{+}1)$-truncated reduced homotopy
\begin{equation*}
    \tilde{r}^i_{\blt}\colon \cosk_n(X_{\blt})_{\blt} \rightarrow \cosk_n(Y_{\blt})_{\blt}
\end{equation*}
between $\cosk_n(f)$ and $\cosk_n(g)$.
\end{srprp}

Moreover, any reduced homotopy induces a unique homotopy.

\begin{srlem}[{\cite[Tag 019L]{sta}}]\label{lem:red-ho-ho-bij}
Let $f,g\colon X_{\blt} \rto Y_{\blt}$ be morphisms of simplicial objects in a category that has finite coproducts.
There is a natural bijection between reduced homotopies and homotopies from $f$ to $g$.
\end{srlem}

We are now equipped to finish the proof of Theorem~(\ref{thm:hcov-sawc-ref}).

\begin{proof}[Proof of Theorem~(\ref{thm:hcov-sawc-ref})]
Let $W_{\blt}$ be a split affine wc hypercovering and $U_{\blt}$ a hypercovering of the pro-\'etale site of a qc scheme.
In Lemma~(\ref{lem:hcov-sawc-ref}) we showed the existence of a morphism $W_{\blt} \rto U_{\blt}$.
Concerning its uniqueness, we are left with showing that any two endomorphisms $f$ and $g$ of $W_{\blt}$ are homotopic.
We inductively construct a reduced homotopy from $f$ to $g$.
\par
In degree $0$ we define $r^0_0 \coloneqq f_0$ and $r^1_0 \coloneqq g_0$.
Suppose that an $n$-truncated reduced homotopy~$r^i_k\colon W_k \rto W_k$ has been constructed.
By Proposition~(\ref{prp:n-trunc-red-ho-cosk}) we obtain an $(n{+}1)$-truncated reduced homotopy $\tilde{r}^i_k\colon
\cosk_n(W_{\blt})_{k} \rto \cosk_n(W_{\blt})_k$ of the $n$-coskeleton that agrees with the $r^i_k$ for $k \leqslant n$.
We use this to extend the reduced homotopy of $W_{\blt}$ to degree $n{+}1$.
First, we define $r^0_{n+1} \coloneqq f_{n+1}$ and $r^{n+2}_{n+1} \coloneqq g_{n+1}$.
Given an index $0 < i \leqslant n{+}1$, we force $r^i_{n+1}$ to be compatible with the degeneracy morphisms by using that
the hypercovering $W_{\blt}$ is split.
On the summand corresponding to the degeneracy $s^j$ we define
\begin{equation*}
    r^i_{n+1} \coloneqq
        \begin{cases*}
            s^j \circ r^i_n & for $j \leqslant i$\\
            s^j \circ r^{i-1}_n & for $j < i$.
        \end{cases*}
\end{equation*}
Lastly, we have to define $r^i_{n+1}$ on the non-degenerate summand $\NDeg W_{n+1}$ and verify the identities involving the
face morphisms.
Since $\NDeg W_{n+1}$ is affine wc, using Corollary~(\ref{cor:char-awc-proj}) we obtain a dashed lift
\begin{equation*}
\begin{tikzpicture}
    \node (NWn+1)  at (0,   0) {$\NDeg W_{n+1}$};
    \node (cWn+11) at (3,   0) {$\cosk_n(W_{\blt})_{n+1}$};
    \node (cWn+12) at (6.5,   0) {$\cosk_n(W_{\blt})_{n+1}\hsmash{,}$};
    \node (Wn+1)   at (6.5, 1.5) {$W_{n+1}$};
    \path[-{Straight Barb[length=2pt]}, font=\scriptsize, line width=0.2mm]
        (NWn+1)  edge (cWn+11)
        (cWn+11) edge node[below] {$\tilde{r}^i_{n+1}$} (cWn+12)
        (Wn+1)   edge (cWn+12);
    \path[dashed, -{Straight Barb[length=2pt]}, font=\scriptsize, line width=0.2mm]
        (NWn+1)  edge (Wn+1);
\end{tikzpicture}
\end{equation*}
which we denote by $r^i_{n+1}$.
The identities involving the face morphisms follow from the construction of $\tilde{r}^i_{n+1}$ in \cite[10]{bkr19} and
the definition of $r^i_{n+1}$ as its lift.
\par
Thus we obtain a reduced homotopy from $f$ to $g$, which corresponds to a homotopy from $f$ to~$g$ by
Lemma~(\ref{lem:red-ho-ho-bij}).
\end{proof}

\chapter{Cohomology of the pro-\'etale homotopy type}\label{sec:coho-peht}

\begin{srsum}
We define a notion of cohomology for the pro-\'etale homotopy type for a certain class of pro-\'etale sheaves and then show
that it computes their pro-\'etale cohomology.
\end{srsum}

The central tool in this section is the Verdier hypercovering theorem, which allows us to compute the cohomology of
abelian sheaves using the homotopy category of hypercoverings.

\begin{srthm}[{\cite[Th\'eor\`eme V.7.4.1.2]{gro72b}}]\label{thm:verdier-hcov}
Let $\sC$ be a site that admits finite limits and finite coproducts.
For any abelian sheaf $\sF$ on $\sC$ there is a natural isomorphism
\begin{equation*}
    \Homol^p(\sC, \sF) = \colim_{U_{\blt}\in\Ho(\Hy)^{\op}} \Homol^p(\sF(U_{\blt}))
\end{equation*}
for all $p \geqslant 0$.
\end{srthm}

Recall that the \'etale homotopy type computes the \'etale cohomology with values in constant abelian sheaves.

\begin{srrem}
Let $S$ be a locally Noetherian scheme and $G$ an abelian group.
We have a natural isomorphism
\begin{equation*}
    \Homol^p(S_{\et}, \ul{G}) = \colim_{U_{\blt}\in\Ho(\Hy)^{\op}} \Homol^p(\comp(U_{\blt}), G) \eqqcolon \Homol^p(\Pi_{
        \et}(S), G).
\end{equation*}
by \cite[Corollary 9.3]{am69}.
Here the cohomology groups in the middle are those of the simplicial set $\comp(U_{\blt})$.
\end{srrem}

We generalise this for the pro-\'etale site of a qcqs scheme using the subcategory of sheaves characterised in
\cite[Lemma 4.2.13]{bs14}, which we give a name.

\begin{srdef}\label{def:comp-sh}
For any qcqs scheme $S$, a \emph{$\comp$-sheaf} is a sheaf $\sF$ on the site $S_{\pe}$ such that for any morphism $X
\rightarrow Y$ of qcqs objects in $S_{\pe}$ that induces an isomorphism between the spaces of components, the morphism
\begin{equation*}
    \sF(Y) \rightarrow \sF(X)
\end{equation*}
is an isomorphism.
The $\comp$-sheaves form a full subcategory $\comp$-$\Sh(S_{\pe})$ of the category of sheaves.
For an \emph{abelian $\comp$-sheaf} the morphism $\sF(Y) \rightarrow \sF(X)$ is required to be an isomorphism of abelian
groups and the category is denoted by $\comp$-$\AbSh(S_{\pe})$.
\end{srdef}

We have already seen examples of $\comp$-sheaves.

\begin{srexa}
The constant sheaf associated to a totally disconnected space $T$ on the pro-\'etale site of a qcqs scheme is a
$\comp$-sheaf.
In fact, we have $\ul{T}(X) = \ul{T}(\comp(X))$ for any object $X$ in the pro-\'etale site by
Proposition~(\ref{prp:comp-incl-adj}).
Similarly we see that a constant sheaf associated to an abelian totally disconnected group is an abelian $\comp$-sheaf.
\end{srexa}

By an abstract argument the equivalence of \cite[Lemma 4.2.13]{bs14} of $\comp$-sheaves with sheaves of the
pro-\'etale site on the space of components factors through the categories of abelian sheaves.

\begin{srcor}\label{cor:comp-absh-eq}
Assuming $S$ is a qcqs scheme, the functors $\alpha_{\sh}$ and $\alpha^{\sh}$ induced by the morphism of sites from
Proposition~(\ref{prp:con-pe-covs}) give rise to an equivalence
\begin{equation*}
    \alpha_{\sh} \simeq \alpha^{\sh}\colon \AbSh(\comp(S)_{\pe}) \rightleftarrows \comp{\text-}\AbSh(S_{\pe})
\end{equation*}
of categories by restriction of $\alpha^{\sh}$.
Moreover, for any abelian $\comp$-sheaf $\sF$ and qcqs object $X$ in the pro-\'etale site there is a natural isomorphism
\begin{equation*}
    \alpha^{\sh}\sF(\comp(X)) = \sF(X).
\end{equation*}
\end{srcor}

\begin{proof}
For any morphism $\varphi\colon \sC \rightarrow \sD$ of sites the induced functors $\varphi_{\sh}$ and $\varphi^{\sh}$ on
the categories of sheaves factor through the respective categories of abelian sheaves.
On the level of sheaves of sets the equivalence is proved in \cite[Lemma 4.2.13]{bs14}.
\par
We conclude the second claim from Lemma~(\ref{lem:comp-pe-covs}) and the defining property of abelian $\comp$-sheaves.
\end{proof}

With Corollary~(\ref{cor:comp-absh-eq}) in place, we naturally extend the definition of the cohomology of the
pro-\'etale homotopy type to $\comp$-sheaves.

\begin{srdef}\label{def:coho-peht-comp-sh}
Given a qcqs scheme $S$ and an abelian $\comp$-sheaf $\sF$, the \emph{cohomology of the pro-\'etale homotopy type} is
defined by
\begin{equation*}
    \Homol^p(\Pi_{\pe}(S), \sF) \coloneqq \colim_{W_{\blt}\in\Ho(\Hyawc)^{\op}} \Homol^p(\alpha^{\sh}\sF(\Pi_{\pe}(S)
        (W_{\blt}))).
\end{equation*}
In fact, by Theorem~(\ref{thm:sawc-hcov-ini}) it is enough to consider a single split affine wc hypercovering.
\end{srdef}

We verify that this generalised definition of the cohomology of the pro-\'etale homotopy type computes the pro-\'etale
cohomology.

\begin{srthm}\label{thm:coho-peht-comp-sh-pe-coho}
For any qcqs scheme $S$ and abelian $\comp$-sheaf $\sF$ there is a natural isomorphism
\begin{equation*}
    \Homol^p(\Pi_{\pe}(S), \sF) = \Homol^p(S_{\pe}, \sF).
\end{equation*}
\end{srthm}

\begin{proof}
By unpacking the definitions and applying Corollary~(\ref{cor:comp-absh-eq}) we have
\begin{equation*}
    \Homol^p(\alpha^{\sh}\sF(\Pi_{\pe}(S)(W_{\blt}))) = \Homol^p(\alpha^{\sh}\sF(\comp(W_{\blt}))) = \Homol^p(\sF(W_{\blt}))
\end{equation*}
for any affine wc hypercovering $W_{\blt}$.
After taking the colimit over all such hypercoverings we obtain the cohomology of the pro-\'etale homotopy type of $S$.
Using Corollary~(\ref{cor:wc-hcovs-ini}) we conclude
\begin{equation*}
    \colim_{W_{\blt}\in\Ho(\Hyawc)^{\op}} \Homol^p(\sF(W_{\blt})) = \colim_{U_{\blt}\in\Ho(\Hy)^{\op}} \Homol^p(\sF(U_{
        \blt})).
\end{equation*}
This equals the pro-\'etale cohomology of $\mathscr{F}$ by Theorem~(\ref{thm:verdier-hcov}).
\end{proof}

\chapter{Pro-\'etale homotopy type of a field}\label{sec:peht-field}

\begin{srsum}
We give some ideas on how to describe the pro-\'etale homotopy type of a field.
Although we do not succeed in full generality, we compute the pro-\'etale homotopy type of the real numbers.
\end{srsum}

\begin{srrem}
Let $k$ be a field and $\sep{k}$ a separable closure in an algebraic closure.
Contrary to the \'etale situation the natural covering $\Spec(\sep{k}) \rto \Spec(k)$ is an object in the pro-\'etale site.
Since the pro-\'etale fundamental group of $k$ is given by $\Gal(\sep{k}/k)$, we expect the pro-\'etale homotopy type to
be its classifying space.
\end{srrem}

We first show that the separable closure is an affine wc covering.

\begin{srlem}\label{lem:k-sep-awc}
Given a field $k$ and a separable closure $\sep{k}$ in an algebraic closure, the affine scheme $\Spec(\sep{k})$ is a wc
object in the pro-\'etale site of $k$.
\end{srlem}

\begin{proof}
Let $\Spec(R) \rightarrow \Spec(\sep{k})$ be an affine covering in the pro-\'etale site of $k$.
Choose a point~$p$ of $\Spec(R)$.
By \cite[Tag 092R]{sta} its residue field $\kappa(p)$ is a separable extension of $k$.
Hence there is a $k$-morphism $\kappa(p) \rightarrow \sep{k}$.
The composite $\sep{k} \rightarrow R \rightarrow \kappa(p) \rightarrow \sep{k}$ is a $k$-automorphism of $\sep{k}$ and thus
has an inverse.
By composing the morphism $R \rightarrow \sep{k}$ with this inverse we obtain a section of the morphism $\sep{k} \rightarrow
R$.
We conclude that $\Spec(\sep{k})$ is wc.
\end{proof}

\begin{srrem}
Let $k$ be a field and $\sep{k}$ a separable closure in an algebraic closure.
Our goal is to use the affine wc covering by $\sep{k}$ to construct an affine wc hypercovering of $k$.
We first examine the natural hypercovering $\cosk_0(\Spec(\sep{k}))_{\blt}$.
\end{srrem}

We can describe the underlying space of the $0$-coskeleton of a separable closure.

\begin{srthm}\label{thm:comp-spc-cosk0-k-sep}
Given a field $k$ and a separable closure $\sep{k}$ in an algebraic closure, the space underlying $\cosk_0(\Spec(\sep{k}))_n$
is naturally isomorphic to the $n$-fold product of $\Gal(\sep{k}/k)$.
\end{srthm}

\begin{proof}
We follow a suggestion by Annette Huber-Klawitter.
Since the field $\sep{k}$ is a filtered colimit of the finite Galois-extensions of $k$ that are contained in $\sep{k}$,
we may write
\begin{equation*}
    \Spec(\sep{k} \otimes_k \sep{k}) = \lim_{l/k} \Spec(\sep{k} \otimes_k l) = \lim_{l/k} \coprod_{\Gal(l/k)} \Spec(\sep{k}).
\end{equation*}
Any connected subset of $\Spec(\sep{k} \otimes_k \sep{k})$ gets mapped to a point under the projections of the filtered limit
since the involved spaces are discrete.
Any ideal $\mfi$ of $\sep{k} \otimes_k \sep{k}$ can be written as
\begin{equation*}
    \mfi = \colim_{l/k} (\mfi \cap (\sep{k} \otimes_k l)).
\end{equation*}
In particular, this holds for any prime ideal and we conclude that any connected subset of the scheme $\Spec(\sep{k}
\otimes_k \sep{k})$ consists of a single point.
Furthermore, we have $\mfi$ is contained in $\mfp$ if and only if $\mfi \cap (\sep{k} \otimes_k l)$ is contained in $\mfp
\cap (\sep{k} \otimes_k l)$ for all finite Galois-extensions $l/k$ in $\sep{k}$.
Hence the space underlying the spectrum of $\sep{k} \otimes_k \sep{k}$ is naturally isomorphic to $\Gal(\sep{k}/k)$.
These calculations work analogously for an $(n{+}1)$-fold tensor product.
\end{proof}

At first sight the situation seems promising.

\begin{srrem}
Let $k$ be a field and $\sep{k}$ a separable closure in an algebraic closure.
After verifying the simplicial structures agree the space of components of $\cosk_0(\Spec(\sep{k}))_{\blt}$ should be given
by the classifying space of $\Gal(\sep{k}/k)$.
We run into the problem that $\cosk_0(\Spec(\sep{k}))_{\blt}$ need not be wc.
By \cite[Theorem 1.8]{bs14} a scheme is affine wc if and only if its space of components is w-local, its local rings at closed
points are strictly Henselian and its space of components is extremally disconnected.
The first two properties are satisfied in each degree since a profinite set is w-local by \cite[Example 2.1.2]{bs14} and
the local rings are given by $\sep{k}$.
The problem lies in the space of components.
\end{srrem}

We give an example that illustrates the problem.

\begin{srexa}
Let $k$ be a finite field and $\sep{k}$ a separable closure.
By \cite[Examples 1.3.7.1]{sza09} the Galois-group of $\sep{k}/k$ is $\hat{\ZZ}$, the profinite completion of the integers,
which is not extremally disconnected.
Assume the contrary and note that we can write $\hat{\ZZ}$ as a product over all $p$-adic integers.
Since the projections of a product of spaces are open, we conclude that the $p$-adic integers are extremally disconnected.
This is absurd; see \cite{mo328070}.
\end{srexa}

\begin{srrem}
Even if the Galois-group were extremally disconnected, it is known that a product of two non-discrete extremally disconnected
groups is not extremally disconnected.
See for example \cite[Warning 2.6]{sch19} or \cite[Problem 4.5.A]{at08}.
\end{srrem}

For discrete groups we do not run into these problems.

\begin{srexa}\label{exa:peht-reals}
The Galois-group of the complex numbers over the real numbers is $\ZZ/2\ZZ$ and thereby discrete.
Hence the $0$-coskeleton of the natural covering by the complex numbers is affine wc.
Furthermore, in this situation the component functor is the one from the \'etale setting and we conclude that
Theorem~(\ref{thm:comp-spc-cosk0-k-sep}) extends to an isomorphism of simplicial spaces.
By Theorem~(\ref{thm:sawc-hcov-ini}) the pro-\'etale homotopy type of the real numbers is given by the classifying
space of $\ZZ/2\ZZ$.
\end{srexa}

We have the hope that Lemma~(\ref{lem:con-w-c}) can be used to construct an affine wc hypercovering in the pro-\'etale site
of a field.

\chapter{Refined pro-\'etale homotopy type}\label{sec:refd-peht}

\begin{srsum}
In the following we define the pro-\'etale homotopy type as a pro-object of condensed sets using the machinery introduced
by Ilan Barnea and Tomer M.~Schlank in \cite{bs16}.
\end{srsum}

First, we recall the notion of condensed set.

\begin{srrem}
A condensed set is a sheaf on the pro-\'etale site of an algebraically closed field and we denote their category by $\Cond$.
Equivalently it is a sheaf on the site of extremally disconnected spaces with finite surjective families as coverings.
The category of condensed sets provides a natural setting to combine topology with homological algebra.
Indeed, there is a faithful functor $\ul{\earg}\colon \Topo \rto \Cond$ that becomes full when restricted to compactly
generated spaces; see Lemma~(\ref{lem:sheaf-assoc-to-space}) and \cite[Proposition 1.7]{sch19}.
\end{srrem}

We review the definition of the refined \'etale homotopy type from \cite{bs16}.

\begin{srrem}
In \cite{bs16} Ilan Barnea and Tomer M.~Schlank introduce the notion of weak fibration category and show that any category
of simplicial sheaves admits such a structure; see~\cite[Section 7.1]{bs16}.
The weak equivalences and fibrations in this structure are the combinatorial weak equivalences and local fibrations considered
in \cite{jar87} by John F.~Jardine.
This weak fibration category induces a model structure on the category of pro-objects.
The authors use this new model structure to define a refinement of the \'etale homotopy type from \cite{am69}.
\par
Let $\mfS$ be a topos such that the left adjoint $\Gamma^*$ of the natural geometric morphism $\Gamma\colon \Set \rlto \mfS$
admits a left adjoint $\comp_{\mfS}$.
By \cite[Proposition 8.1]{bs16} the induced functor $\Gamma^*$ on simplicial objects is a weak right Quillen functor with
respect to the weak fibration structures on the categories of simplicial objects in $\Set$ and $\mfS$.
The induced functors $\Pro(\comp_{\mfS})$ and $\Pro(\Gamma^*)$ on pro-objects form a Quillen
adjunction
\begin{equation*}
    \Pro(\comp_{\mfS}) \dashv \Pro(\Gamma^*)\colon \Pro(\simp\mfS) \rlto \Pro(\sSet)
\end{equation*}
by \cite[Corollary 5.4]{bs16}.
The authors now define the object $\LDer\Pro(\comp_{\mfS})(*)$ in the homotopy category of $\Pro(\sSet)$ to be the refined
homotopy type of $\mfS$.
Here $\LDer\Pro(\comp_{\mfS})$ denotes the left derived functor induced by the Quillen adjunction and the homotopy category
is the one induced by the model structure.
In \cite[Proposition 8.4]{bs16} the authors compare their definition to the one from \cite{am69} for the \'etale
topos of a locally Noetherian scheme.
\end{srrem}

With some modifications we can apply the theory from \cite{bs16} to the pro-\'etale topos of a qcqs scheme.

\begin{srrem}
Let $S$ be a qcqs scheme.
We consider the full subcategory $S_{\awc, \pe}$ of affine wc objects in the pro-\'etale site of $S$, which carries a natural
structure of site.
The space of components induces a functor $\comp\colon S_{\awc, \pe} \rto \Extr$, which respects coverings but need not
be a morphism of sites.
Still, we have a pair $\comp_{\psh} \dashv \comp^{\psh}$ of adjoint functors between the categories of presheaves using
\cite[Theorem I.2.3.1]{tam94}.
Note that the inclusion of $S_{\awc, \pe}$ into the pro-\'etale site of $S$ induces an equivalence on the categories of sheaves.
Using the sheafification functor we get a pair
\begin{equation*}
    \comp_{\sh} \dashv \comp^{\sh}\colon \Sh(S_{\pe}) \rlto \Cond
\end{equation*}
of adjoint functors as in \cite[64]{tam94}.
Note that this need not be a geometric morphism due to the defect mentioned above.
The left adjoint should be thought of as an extension of the component functor to the pro-\'etale topos.
\end{srrem}

In fact, we show that the left adjoint is given by the component functor on affine wc objects.
We thank Timo Richarz for pointing us towards this fact.

\begin{srlem}\label{lem:comp-sh-repr}
Let $S$ be a qcqs scheme.
Given an affine wc object $W$ in the pro-\'etale site of $S$, there is a natural isomorphism $\comp_{\sh}\rep_W =
\ul{\comp(W)}$.
\end{srlem}

\begin{proof}
Let $A$ be a condensed set.
By Lemma~(\ref{lem:shfctn-awc-objs}) we have $\comp^{\sh}A(W) = A(\comp(W))$.
We conclude that there are natural bijections
\begin{equation*}
    \Hom(\comp_{\sh}\rep_W, A) = \Hom(\rep_W, \comp^{\sh}A) = \comp^{\sh}A(W) = \Hom(\ul{\comp(W)}, A),
\end{equation*}
which proves the claim.
\end{proof}

We need to verify manually that the right adjoint preserves epimorphisms to conclude that it is a weak right Quillen functor.

\begin{srlem}\label{lem:comp-up-s-preserv-epis}
For any qcqs scheme $S$ the functor $\comp^{\sh}$ preserves epimorphisms.
\end{srlem}

\begin{proof}
Let $A \rto B$ be an epimorphism of condensed sets.
Since the pro-\'etale site has enough affine wc objects, a morphism of sheaves is an epimorphism if and only if its evaluation
at any such object is surjective.
By Lemma~(\ref{lem:shfctn-awc-objs}) applying the functor $\comp^{\sh}$ and then evaluating at an affine wc object $W$ gives
back the morphism $A(\comp(W)) \rto B(\comp(W))$.
The claim follows since the space of components preserves affine wc objects by Proposition~(\ref{prp:awc-pe-site-k-alg}).
\end{proof}

We are now ready to prove the central theorem of this section.

\begin{srthm}\label{thm:comp-quil-adj}
Given a qcqs scheme $S$, the functors $\comp_{\sh}$ and $\comp^{\sh}$ induce a pair
\begin{equation*}
    \Pro(\comp_{\sh}) \dashv \Pro(\comp^{\sh})\colon \Pro(\simp\Sh(S_{\pe})) \rlto \Pro(\sCond)
\end{equation*}
of Quillen-adjoint functors with respect to the model structure induced by the natural weak fibration categories on the
pro-\'etale topoi; see \cite[Section 7.1]{bs16}.
\end{srthm}

\begin{proof}
The functor $\comp^{\sh}$ preserves limits as a right adjoint and epimorphisms of sheaves by
Lemma~(\ref{lem:comp-up-s-preserv-epis}).
We conclude by~\cite[Lemma 7.8]{bs16} and \cite[Lemma 7.9]{bs16} that it preserves local (trivial) fibrations.
Hence it is a weak right Quillen functor.
We see that $\Pro(\comp_{\sh})$ is left adjoint to $\Pro(\comp^{\sh})$ using \cite[Lemma 5.8.2]{bs16}.
By \cite[Corollary 5.4]{bs16} we have the claimed Quillen adjunction.
\end{proof}

We give a refined definition of the pro-\'etale homotopy type and compare it with our earlier definition.

\begin{srdef}
For any qcqs scheme $S$ its \emph{refined pro-\'etale homotopy type} is defined to be the object
\begin{equation*}
    \tilde{\Pi}_{\pe}(S) \coloneqq \LDer\Pro(\comp_{\sh})(*)
\end{equation*}
in the homotopy category of $\Pro(\sCond)$.
Here $\LDer\Pro(\comp_{\sh})$ denotes the left derived functor induced by the Quillen adjunction from
Theorem~(\ref{thm:comp-quil-adj}) and the homotopy category is the one obtained by inverting weak equivalences.
\end{srdef}

Since we are working with simplicial condensed sets instead of simplicial sets, we have to transfer some results to our
setting.

\begin{srrem}\label{rem:ho-eq-sprof-vs-scond}
We write $\sCond$ for the natural weak fibration category of simplicial condensed sets and denote its homotopy category
obtained by inverting weak equivalences by $\How(\sCond)$; see \cite[Section 7]{bs16}.
By $\Pro(\sCond)$ we denote the model category of pro-objects induced by the weak fibration structure; see \cite[Section 4]{bs16}.
\par
For any fibrant object $A_{\blt}$ in $\sCond$ we have the natural functorial path object
\begin{equation*}
    A_{\blt} = A_{\blt}^{\Delta[0]} \rto A_{\blt}^{\Delta[1]} \rto A_{\blt}^{\Delta[0] \bcoprod \Delta[0]} = A_{\blt} \times
        A_{\blt}
\end{equation*}
induced by the simplicial structure; see proof of \cite[Proposition 8.4]{bs16}.
The exponentials come from the presheaf exponential $\sF^{\sP}(U) \coloneqq \Hom(\rep_U \times \sP, \sF)$, which is a sheaf
if $\sF$ is one.
The notion of homotopy induced by this functorial path object is the same as the one induced by the simplicial condensed
set $A_{\blt} \otimes \Delta[1]$ defined in \cite[837]{bs16}.
It is the sheafification of the simplicial presheaf $E \mapsto A_{\blt}(E) \times \Delta[1]$ on the category of extremally
disconnected spaces.
Therefore we have
\begin{equation*}
    A_{\blt} \otimes \Delta[1] = A_{\blt} \times \ul{\Delta[1]}.
\end{equation*}
Any homotopy equivalence in $\sCond$ is a weak equivalence.
The functor $\ul{\earg}\colon \sProf \rto \sCond$ sending a simplicial profinite set to its constant simplicial condensed
set preserves binary products and thus homotopies.
All in all, we obtain a functor $\ul{\earg}\colon \Ho(\sProf) \rto \How(\sCond)$.
\end{srrem}

\begin{srrem}
We note that a weak equivalence in $\Pro(\sCond)$ is a morphism that is isomorphic to a natural transformation that
is level-wise a weak equivalence in $\sCond$.
In particular, the natural functor on the categories of pro-objects induced by $\upgamma\colon \sCond \rto \How(\sCond)$
sends weak equivalences to isomorphisms.
\end{srrem}

We are now able to follow the proof of \cite[Proposition 8.4]{bs16} almost verbatim.

\begin{srthm}\label{thm:peht-vs-refd-peht}
Given a qcqs scheme $S$, there is a natural isomorphism
\begin{equation*}
    \upgamma \circ \tilde{\Pi}_{\pe}(S) = \ul{\Pi_{\pe}(S)}
\end{equation*}
between the refined pro-\'etale homotopy type and the pro-\'etale homotopy type considered as a pro-object in $\How(\sCond)$.
Here $\upgamma$ denotes the natural functor $\sCond \rto \How(\sCond)$.
\end{srthm}

\begin{proof}
We denote $\mfS \coloneqq \Sh(S_{\pe})$ and note that there is a natural functorial path object for $\simp\mfS$ by
Remark~(\ref{rem:ho-eq-sprof-vs-scond}).
By applying the discussion in \cite[828]{bs16} to the final object $*$ in $\Pro(\simp\mfS)$ we obtain a cofibrant replacement
\begin{equation*}
    j\colon \sA_* \rto \widehat{\simp\mfS}_{\fw} \subset \simp\mfS
\end{equation*}
Here $\simp\mfS_{\fw}$ denotes the full subcategory of $\simp\mfS$ consisting of the fibrant and contractible objects.
The category $\widehat{\simp\mfS}_{\fw}$ has the same objects as $\simp\mfS_{\fw}$, but only morphisms that are fibrations
and weak equivalences; see \cite[832]{bs16} for more details.
\par
We can use the cofibrant replacement of $*$ to compute the refined pro-\'etale homotopy type: it is given by $\comp_{\sh}
\circ j$.
Moreover, we have a commutative diagram
\begin{equation*}
\begin{tikzpicture}
    \node (A) at (0, 1.5) {$\sA_*$};
    \node (sXhfw) at (3, 1.5) {$\widehat{\simp\mfS}_{\fw}$};
    \node (pisXhfw) at (3, 0) {$\comp_0(\widehat{\simp\mfS}_{\fw})$};
    \node (sXfw) at (6, 1.5) {$\simp\mfS_{\fw}$};
    \node (pisXfw) at (6, 0) {$\comp_0(\simp\mfS_{\fw})$};
    \node (sCond) at (9, 1.5) {$\sCond$};
    \node (HosCond) at (9, 0) {$\How(\sCond)\hsmash{.}$};
    \path[-{Straight Barb[length=2pt]}, font=\scriptsize, line width=0.2mm]
        (A)       edge node[above] {$j$} (sXhfw)
        (sXhfw)   edge (sXfw)
                  edge (pisXhfw)
        (pisXhfw) edge (pisXfw)
        (sXfw)    edge node[above] {$\comp_{\sh}$} (sCond)
                  edge (pisXfw)
        (pisXfw)  edge node[below] {$\comp_{\sh}$} (HosCond)
        (sCond)   edge node[right] {$\upgamma$} (HosCond);
\end{tikzpicture}
\end{equation*}
as in the proof of \cite[Proposition 8.4]{bs16}.
Here $\comp_0(\simp\mfS_{\fw})$ has the same objects as $\simp\mfS_{\fw}$, but homotopic morphisms are identified.
The same is true for $\comp_0(\widehat{\simp\mfS}_{\fw})$; see \cite[Definition 6.15]{bs16}.
Note that we use the notation ``$\comp_0$'' to avoid confusion with the space of components.
\par
By \cite[833]{bs16} the composite $\sA_* \rto \comp_0(\simp\mfS_{\fw})$ is initial and thus $\upgamma \circ \comp_{\sh}
\circ j$ is naturally isomorphic to the pro-object $\comp_0(\simp\mfS_{\fw}) \rto \How(\sCond)$.
To obtain $\Pi_{\pe}(S)$ considered as a pro-object in $\How(\sCond)$, we need to restrict to the homotopy category of
hypercoverings.
The inclusion $\Ho(\Hy(S_{\pe})) \rto \comp_0(\simp\mfS_{\fw})$ is initial by \cite[Lemma 2.2]{jar94}.
Here we consider a hypercovering as a constant simplicial sheaf by virtue of the associated representable sheaf.
\par
We have proved that $\upgamma \circ \tilde{\Pi}_{\pe}(S)$ is naturally isomorphic to $\Ho(\Hy(S_{\pe})) \rto \How(\sCond)$
as a pro-object.
By Corollary~(\ref{cor:wc-hcovs-ini}) we may further restrict to affine wc hypercoverings.
On those the functor $\comp_{\sh}$ is given by the space of components followed by $\ul{\earg}$ by Lemma~(\ref{lem:comp-sh-repr}).
\end{proof}

We also obtain a natural definition for the homotopy groups of the pro-\'etale homotopy type.

\begin{srrem}\label{rem:pro-prof-pi-n}
For a pointed qcqs scheme our first version of the pro-\'etale homotopy type is profinite.
Thus we may use the theory in \cite{qui08} by Gereon Quick to define pro-profinite homotopy groups associated to the pro-\'etale
homotopy type.
But as noted in \cite[70]{bs14}, the pro-\'etale fundamental group need not even be pro-discrete.
\end{srrem}

\begin{srrem}
Given a simplicial presheaf $\sF_{\blt}$, we obtain a notion of homotopy groups fibred over $\sF_0$ as explained in
\cite[833-834]{bs16}.
For any simplicial set $E_{\blt}$, we consider
\begin{equation*}
    \comp_n(E_{\blt}) \coloneqq \coprod_{e\in E_0} \comp_n(E_{\blt}, e),
\end{equation*}
where we take the coproduct over the usual pointed homotopy groups of simplicial sets.
We have a natural map $\comp_n(E_{\blt}) \rto E_0$ that is functorial in $E_{\blt}$.
For any simplicial condensed set $A_{\blt}$ we obtain a presheaf $\comp_n \circ A_{\blt}$ of groups fibred over $A_0$ on
the site of extremally disconnected spaces.
By this we mean that for any extremally disconnected space $P$, the set $\comp_n(A_{\blt}(P))$ is a group fibred over the
set $A_0(P)$.
Note that contrary to the general case this is already a sheaf, and thereby a condensed set, since the functor $\comp_n$
preserves binary products.
We think of this as the $n$th condensed homotopy group of $A_{\blt}$ fibred over $A_0$.
\par
Given a qcqs scheme $S$, we obtain fibred homotopy pro-groups $\comp_n \circ \tilde{\Pi}_{\pe}(S)$ associated to the refined
pro-\'etale homotopy type.
By definition weak equivalences of simplicial sheaves induce isomorphisms on the fibred homotopy groups.
Hence we obtain a natural isomorphism
\begin{equation*}
    \comp_n \circ \tilde{\Pi}_{\pe}(S) = \comp_n \circ \ul{\Pi_{\pe}(S)}
\end{equation*}
by Theorem~(\ref{thm:peht-vs-refd-peht}).
\end{srrem}

We unravel the definitions in the pointed setting.

\begin{srrem}\label{rem:pi-n-ptd-scond}
Given a simplicial condensed set $A_{\blt}$ and a point $a$ in $A_0(*)$, we can describe the fibred homotopy groups of
$A_{\blt}$ above $a$.
They are given as the sheaves $E \mapsto \comp_n(A_{\blt}(E), a_E)$ on the site of extremally disconnected spaces, where
$a_E$ denotes the image of $a$ under the natural map $A_0(*) \rto A_0(E)$.
In particular, any pointed simplicial profinite set induces a pointed simplicial condensed set in the sense above.
Moreover, this construction is functorial in $A_{\blt}$ and the point $a$.
\end{srrem}

This motivates our definition of the pro-\'etale homotopy groups of a pointed qcqs scheme.

\begin{srdef}
For a pointed qcqs scheme $(S, s)$ the \emph{$n$th pro-\'etale homotopy group} is the pro-object
\begin{equation*}
    \comp_{\pe, n}(S,s) \coloneqq \comp_n \circ \ul{\Pi_{\pe}(S, s)}
\end{equation*}
in the category of condensed groups.
Here the functor $\comp_n$ on the right denotes the condensed homotopy group from Remark~(\ref{rem:pi-n-ptd-scond}).
\end{srdef}

\begin{srrem}
Our definition raises the question of its relation to the pro-profinite homotopy groups from Remark~(\ref{rem:pro-prof-pi-n})
and, more importantly, to the pro-\'etale fundamental group.
At this point in time we have not been able to answer these questions.
\end{srrem}

\printbibliography

\end{document}